# A comment on the combination of the implicit function theorem and the Morse lemma

Emőke Imre

*Óbuda University, Hydro-Bio-Mechanical Systems Research Center, Budapest, Hungary*

**Abstract**

The analytic implicit function theorem is extended. The function *f* of the theorem is integrated with respect to the dependent variable of the implicit function. A geometrical interpretation is given for the sub-geometry of the integral function *F* by using the Morse lemma. The result is used in the analysis of the hierarchical technique related to the minimization inon-linear parameter identification.

## Introduction

**The hierarchical solution of non-linear inverse problems**

The solution of the non-linear inverse problems can be determined by non-linear minimization which is hindered by the fact that the LS merit function has an unimaginable complexity of the *M*-dimensional topography due to the noise, where *M* is the parameter number.

The minima can be *global* or *local*. Two heuristics are used: (i) to find local extrema starting from varying initial values of the independent variables, and then pick the most extreme [1]; or (ii) to perturb a local extremum by taking a step away from it, and then see if the iteration returns to a better point, or to the same one [1]. Another approach is to fit locally a nicer hyper-surface and descending along it [2]. Some additional geometrical difficulties may arise in case of a quasi-degenerated global minimum, which can theoretically be treated by regularization ([3]) needing modified algorithms.

No hierarchical solution method has been reported for the non-linear case despite of the fact that the number of the parameters in the non-linear algorithm can be decreased resulting in less critical points due to the noise and less numerical work besides other advantages (the one-dimensional sections of the merit functions so determined can be used for parameter error estimation and implicit regularization).

**The aim of the paper**

The concept of the hierarchical method [4] in well-known for linear system of equation eg., int he context of linear inverse problems and is started to be introduced into the solution of non-linear inverse problems here. The use of the hierarchical solution is visualized as follows (Fig.1).

We may assume that $F : \mathbf{R}^M = \mathbf{R}^{n+m} \rightarrow \{0 \cup \mathbf{R}^+\}$ is an analytic Least Squares merit function which has a single minimum $p_{\min}$ within the parameter domain. Let us assume that the parameter space is split into direct sum $\boldsymbol{p}=(\boldsymbol{x},\boldsymbol{y})$ and the solution of the inverse problem is split into two, smaller dimensional parts (in Fig. 1: $x^1$, $x^2$), accordingly.

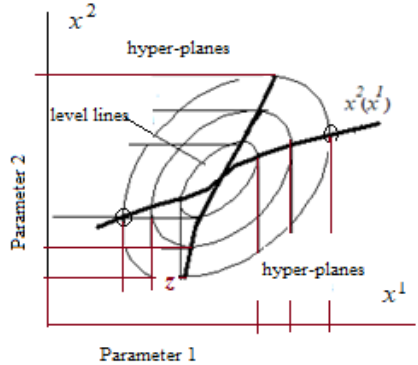

Figure 1. The implicit function $x^2(x^1)$ is defined by $f(x, y) = f(x^1, x^2) = F'_y (x^1, x^2) = 0$. The implicit function $x^2(x^1)$ is the inverse image of the set of the conditional minimum points at fixed $x=x^1$ values. The minimum of $F(x^1, x^2)$ with respect to the two parameters is the same as the minimum of the 'minimal section' $F[x^1, x^2(x^1)]$ if the sub-minima along each vertical ($y=x^2$ directional) hyper-plane is unique.

One kind of minimization happens along the **y** directional parallel plane sections of the parameter space: the $y^2$ solution part is searched in every fixed value of **x** by sub-minimization (using minimization method 1). The set of the so determined sub-minima are found along the graph of an implicit function **y(x)** related to the condition $F'_y = f = 0$. The so called minimal section $F[x,y(x)]$ is minimized (with method 2).

The paper treats the condition of the equivalence of the original and the hierarchical inverse problem solutions, by the analysis of the geometry of the merit function along the **y** planes. It is found that for the equivalence, $F$ is needed to be convex in the **y** direction. The condition is met in case of strictly convex merit functions, or if the model depends linearly on **y** when the merit function is partly strictly convex (ie., in the **y** direction).

Moreover, some consequences and inferences are mentioned. In the strictly convex case, the set of the graphs of the implicit functions related to the various subspaces is isomorphic with the algebraic lattice of the subspaces of the parameter space $\mathbf{R}^M$. The geometry and the analytic nature of the one-dimensional minimal sections can be used in reliability testing and in regularization.

## The extension of the implicit function theorem

### Assumptions, concepts, theorems

The basic concepts used are as follows [5 to 8]. For the smooth function $F : \mathbf{R}^M \to \mathbf{R}$ the points where the first derivative vanishes are called critical points and their images are called critical values. If at a critical point the matrix of second derivative (the Hessian matrix) is non-singular, then it is called a non-degenerate critical point; if the Hessian is singular then it is a degenerate critical point. Some basic theorems are presented which are coupled in the extended implicit function theorem.

Morse lemma

Let $p$ be a non-degenerate critical point of $F$. Then there exists a chart in a neighborhood U of $p$ such that for all $x_i(p) = 0$ where $x_i$ constitute a local map around $p$ and

$$F(x) = F(p) - x_1^2 - \ldots - x_\gamma^2 + x_{\gamma+1}^2 \ldots - x_n^2 \tag{1}$$

throughout U. Here $\gamma$ is equal to the index of $F$ at $p$ (the number of the negative eigenvalues of the Hessian matrix). As a corollary of the Morse lemma, one sees that non-degenerate critical points are isolated.

The last Morse inequality on $S^n$, relating the number of critical points with index $k$ and the Euler characteristic of $S^n$ – being an equality – is as follows:

$$\sum_{k=0}^{n}(-1)^k \text{ (k-index critical points)} = 1 + (-1)^n \tag{2}$$

From this, the equation for $D^n$ (n-dimensional ball) can be derived. A 1 must be subtracted from both sides if a minimum (critical point with index 0) and $(-1)^n$ if a maximum (critical point with index n) is left out. If the function "curls up" at the edge of the ball (the gradient is outward at the edge):

$$\sum_{k=0}^{n}(-1)^k \text{ (k-index critical points)} = 1 \tag{3}$$

In the strict convex case, the edge $D^n$ (n-dimensional ball, which is a half $S^n$) has a single critical point with index of 0.

The implicit function theorem

Let $f: \mathbf{R}^{n+m} \to \mathbf{R}^m$ be an analytic function. $\mathbf{R}^{n+m}$ is the direct sum of $\mathbf{R}^n$ and $\mathbf{R}^m$, a point of this is $(\mathbf{x}, \mathbf{y}) = (x_1, \ldots, x_n, y_1, \ldots, y_m)$. Starting from the given function $f$, the goal is to construct a function $g: \mathbf{R}^n \to \mathbf{R}^m$ whose graph $(\mathbf{x}, g(\mathbf{x}))$ is precisely the set of all $(\mathbf{x}, \mathbf{y})$ such that $f(\mathbf{x}, \mathbf{y}) = 0$.
Fix a point $(\mathbf{a}, \mathbf{b}) = (a_1, \ldots, a_n, b_1, \ldots, b_m)$ with $f(\mathbf{a}, \mathbf{b}) = \mathbf{0}$, where $\mathbf{0}$ is the element of $\mathbf{R}^m$. If the matrix $[(\partial f_i/\partial y_j)(\mathbf{a}, \mathbf{b})]$ is invertible, then there exists an open set $U$ containing $\mathbf{a}$, an open set $V$ containing $\mathbf{b}$, and a unique function $g: U \to V$ such that

$$\{\mathbf{x}, g(\mathbf{x}) | \mathbf{x} \in U\} = \{(\mathbf{x}, \mathbf{y}) | \mathbf{x} \in U, \mathbf{f}(\mathbf{x}, \mathbf{y}) = \mathbf{0}\} \tag{4}$$

Whenever we have the hypothesis that $f$ is analytic inside $U \times V$, then the same holds true for the implicit function $g$ inside $U$.

**Extension**

**Statement 1** *(The extension of the implicit function theorem with geometrical interpretation)*

Let us assume that the parameter space is split into the direct sum $p=(x,y)$ and the conditions of the analytic implicit function theorem are met in $(\mathbf{a},\mathbf{b}) = (a_1, \ldots, a_n, b_1, \ldots, b_m)$, being $f(\mathbf{a}, \mathbf{b}) = \mathbf{0}$, where $\mathbf{0}$ is the element of $\mathbf{R}^m$ and being $f''_y$ regular. As a result, there exists an open set $U$ containing $\mathbf{a}$, and an open set $V$ containing $\mathbf{b}$, and a unique, analytic function $g: U \to V$, $y=g(x)$, $g: \mathbf{R}^n \to \mathbf{R}^m$, such that $f[x,g(x)] = \mathbf{0}$.

(i) Let us assume that the function $f: \mathbf{R}^{n+m} \to \mathbf{R}^m$ in the analytic implicit function theorem has an analytic partial integral function with respect to $y$, denoted by $F: \mathbf{R}^{n+m} \to \mathbf{R}$ such that $f = F'_y$. The $F[x,g(x)]$ is called minimal section of $F$ with respect to $x$. The points $(a,b)$ or $[x,g(x)]$ are non-degenerate, "partial" critical points of $F$ in the $y$ direction, with local sub-geometry determined by the index of $F''_{yy}$. (The points can be the degenerate critical points of $F$).

Explanation:

Let us consider the point $(a,b)$ and it environment U, V of the analytic implicit function theorem. Let us introduce here the partial function $F^x(y) : \mathbf{R}^m \to \mathbf{R}^l$, describing the variation of $y$ at fixed $x$. Its first derivative is equal to $F'_y$ and its second derivative is equal to $F''_{yy}$. The slice function $F^a(y) : \mathbf{R}^m \to \mathbf{R}$ is a Morse function, being its first derivative zero, its Hessian non-degenerate in $b$ having an isolated, non-degenerate critical point there. The same index can be expected along the implicit function due to the nice features of the function $F$.

**Statement 2** *(The application for the hierarchical minimization)*

If $F$ has a non-degenerate minimum in the point $(a,b)$, then for any direct sum decomposition $p=(x,y)$, every possible implicit function locally exists with a local sub-minimum in terms of $y$. However, the implicit function for a special direct sum decomposition $p=(x,y)$ exists globally if $F''_{yy}$ is positive definite globally on the whole parameter domain.

Explanation:

Since $F''_{yy}$ is globally positive definite, the implicit function can locally be extended repetitively (for any $x$). The index of $F^x(y)$ is 0, having a non-degenerate local sub-minimum where $F'_y = f[x,g(x)] = \mathbf{0}$. Moreover, $F$ is strictly convex in the $y$ direction since $F''_{yy}$ is positive definite in the $y$ direction globally. Therefore, $F^x(y)$ has only one critical point, the minimum (the statement also follows from the last Morse inequality).

**Statement 3** *(The case of hierarchical minimization of convex functions on $R^M$)*

The conditions of the implicit function theorem are met in the non-degenerate minimum of $F$ for any direct sum decomposition $p=(x,y)$ and, therefore, the related implicit functions locally exist.

(i) In case of a strictly convex function $F$, the partial derivative $F''_{yy}$ is globally positive definite for any direct sum decomposition $p=(x,y)$. Therefore, every possible minimal section globally and uniquely exists. The set of these is isomorphic with the algebraic lattice of the sub-spaces of $\mathbf{R}^M$ with partial ordering of the containing. If $x'$ is contained by $x$ then the smaller dimensional minimal section of x' is contained by a larger dimensional minimal section of $x$. Therefore, the hierarchical technique can be applied in a repetitive way.

(ii) Let us consider the Euclidean space generated by the graph of $F$. For every z level line there are two coordinate hyper-planes being normal to the coordinate direction of parameter $x^i$ and tangent to the given z level line (Fig. 1). Such hyper-plane pair series bracket decreasing z level lines. The tangent point pairs are the points of the 1-dimensional minimal section of $x^i$, the interval defined by the point pairs $[x^{i\,z-}, x^{i\,z+}]$ is the orthogonal projection of the sub-level set z of F onto the axis of the parameter $x^i$. The intervals decrease with z, in the minimum the two points coincide (Fig. 1).

Corrolary *(properties of the 1-dimensional minimal sections of the parameters)*

The 1-dimensional minimal section $F[[x^i,\mathbf{y}(x^i)]$ is the smallest descent line with respect parameter $x^i$ that proceed in the positive/negative $x^i$ directions toward the global minimizer $p_{min}$. The 1-dimensional minimal section $F[[x^i,\mathbf{y}(x^i)]$ is a 'deepest' sensitivity section (ie., a section taken along a simple curve passing through the global minimum, parametrised by $x^i$).

# Further comments on the hierarchical minimization

## Some comments on the equivalence condition

We have found that if the analytic *F* is partially, globally strictly convex in the *y* direction, then the implicit function of *x* globally and uniquely exists with a sub-minimum of *F* in the *y* direction. The two kinds of minimisation (ie., the *M*-dimensional unconditional and the hierarchical with a *J* dimensional sub-minimisation with respect *y* and the *M-J* dimensional minimisation of the minimal section of *x* with respect *x*) have the same solution.

As a practical significance, it follows that only those parameters are important in the point of view of non-linear minimisation, that influence the model non-linearly. The linear model-part can be 'eliminated'. In the case when only one single parameter influences the model non-linearly, the non-linear minimisation can be solved even by one-dimensional bracketing since all parameters can be eliminated by linear sub-minimisation automatically.

It can be noted that a global implicit function may exist for such partly convex function *F* that has more than one minimum. In case of two distinct minima without common *x* coordinate value, a global implicit function and the related minimal section with respect to *x* will cross both minima.

## Algorithmic steps in a coordinate grid bracketing of convex functions

An implicit function point is a conditional minimum point of *F* and can be determined as follows. Let us consider the orthogonal projection of the graph of *F*(*x*, *y*) onto the *F* - *x* coordinate plan in the total Euclidean space generated by the graph. The critical values of the map - where the fold of the projection is found - are related to the condition $F_{y}'(x, y) = f(x, y) = 0$ which defines the implicit function.

Using this, an algorithm can be elaborated for the determination of the implicit function and the hierarchical minimisation. The unique sub-minimum at the slice function $Fx(y) : R^m \to R$ for fixed x can be determined by projection, minimum search with respect to y and from the minimum, the graph point of the implicit function can be determined by inverse image operation.

A one-dimensional minimum can be bracketed only when there is a *triplet* of points $a < b < c$ (or $c < b < a$), such that $f(b)$ is less than both $f(a)$ and $f(c)$. In this case we know that the function (if it is non-singular) has a minimum in the interval ($a; c$). It is possible to bracket a one-dimensional minimum but there is no analogous procedure in the general, multi-dimensional case [1] unless the function is strictly convex, since hyper-planes generated by coordinate mesh points can enclose decreasing convex level lines.

The global of minimum of strictly convex merit functions can be bracketed by triplets in every coordinate direction, by using a coordinate mesh. Using the same coordinate mesh, by applying the foregoing algorithmic step in each mesh point beyond the global minimum search, the determination of the one-dimensional minimal sections is also possible.

## Regularization

Some additional geometrical difficulties beyond the noise may arise when the merit function is distorted being the minimum quasi-degenerated (eg., the testing time can be "too short" in the measurements in case of the fitting of a model which is the injective solution of a PDE, depending on the time variable). The problem can be treated by regularization ([1, 3]). Explicit regularization is regularization whenever one explicitly adds a term to the optimization problem. These terms could be priors, penalties, or constraints. Explicit regularization is commonly employed with ill-posed optimization problems.

The regularization term, or penalty, imposes a cost on the optimization function to make the optimal solution unique. Implicit regularization is all other forms of regularization. This includes, for example, early stopping, using a robust loss function, and discarding outliers.

The one-dimensional minimal section related to a given parameter has some special geometry features. It is a smallest descent line, it can be used to determine parameter sensitivity / error. Being analytic, it can be used to solve quasi-degenerate minimum problems with implicit regularization as follows.

In this work a new, implicit method is suggested, based on the analytic feature of the one-dimensional implicit functions which may have injective coordinate functions (eg., the $p_{min}$ can be determined by an injective coordinate function of the one-dimensional implicit function related to a given parameter using a known coordinate of this parameter of the parameter vector $p_{min}$). The suggested method does not fail even in such a case when the convexity is slightly destroyed by model perturbation, and there are two degenerated minimuma.

## Acknowledgement

The discussion with András Stipsicz, András Szűcs, Endre Szabó, József Bodnár is greatly acknowledged.